\documentclass[11pt,twoside,reqno]{amsart}
\usepackage{epsfig,amsmath,amssymb,fancyhdr,url,graphicx,latexsym}

\calclayout

\begin{document}

\leftline{ \scriptsize \it  Journal of Prime Research in Mathematics
Vol. {\bf 8}(2012), 61-75 }

\vspace{1.3cm}

\title
{New Recurrence Relationships between Orthogonal Polynomials which Lead to New Lanczos-type Algorithms}

\author{Muhammad Farooq$^1$, Abdellah Salhi$^2$}
\thanks{ {\enskip
  \enskip $^{1}$Department of Mathematics, University of Peshawar, Khyber Pakhtunkhwa, 25120, Pakistan.. Email: mfarooq@upesh.edu.pk\\
$^{2}$Department of Mathematical Sciences, University of Essex, Wivenhoe Park, Colchester, CO4 3SQ, UK. E-mail: as@essex.ac.uk}}
\begin{abstract}
Lanczos methods for solving $\textit{A}\textbf{x}=\textbf{b}$ consist in constructing a sequence of vectors $(\textbf{x}_k), k=1,...$ such that $\textbf{r}_{k}=\textbf{b}-\textit{A}\textbf{x}_{k}=\textit{P}_{k}(\textit{A})\textbf{r}_{0}$,,
where $\textit{P}_{k}$ is the orthogonal polynomial of degree at most $k$ with respect to the linear functional $c$ defined as $c(\xi^i)=(\textbf{y},\textit{A}^i\textbf{r}_{0})$. Let $\textit{P}^{(1)}_{k}$ be the regular monic polynomial of degree $k$ belonging to the family of formal orthogonal polynomials (FOP) with respect to $c^{(1)}$ defined as $c^{(1)}(\xi^{i})=c(\xi^{i+1})$. All Lanczos-type algorithms are characterized by the choice of one or two recurrence relationships, one for $\textit{P}_{k}$ and one for $\textit{P}^{(1)}_{k}$. We shall study some new recurrence relations involving $\textit{P}_{k}$ and $\textit{P}^{(1)}_{k}$ and their possible combination to obtain new Lanczos-type algorithms. We will show that some recurrence relations exist, but cannot be used to derive Lanczos-type algorithms, while others do not exist at all.
\vskip 0.4 true cm
 \noindent
 {\it AMS SUBJECT} : Primary 65F10.\\
\end{abstract}
\maketitle
\bibliographystyle{plain}

\pagestyle{myheadings} \markboth{\centerline {\scriptsize Muhammad Farooq, Abdellah Salhi}}
         {\centerline {\scriptsize  New Recurrence Relationships between Orthogonal Polynomials}}


\bigskip
\bigskip
\medskip

\section{Introduction} \label{sec:intro} Let
In $1950$, C. Lanczos \cite{50:Lanczos} proposed a method for transforming a matrix into a similar tridiagonal matrix. We know, by Cayley-Hamilton theorem, that the computation of the characteristic polynomial of a matrix and the solution of linear equations are equivalent, Lanczos,  \cite{52:Lanczos}, in $1952$ used his method for solving systems of Linear equations.

Since then, several Lanczos-type algorithms have been obtained and among them, the famous conjugate gradient algorithm of Hestenes and Stiefel,  \cite{52:Hest}, when the matrix is Hermitian, and the bi-conjugate gradient algorithm of Fletcher, \cite{76:Fletcher}, for the general case. In the last three decades, Lanczos algorithms have evolved and different variants have been derived, \cite{95:Baheux,98:Bjorck,93:Brezinski,94:Brezinski,92:Brezinski,99:Brezinski,00:Brezinski,02:Brezinski,00:Calvetti,97:Greenbaum,99:Guennouni,06:Gérard,79:Parlett,85:Parlett,87:Saad,87:Van,94:Ye}.

Although Lanczos-type algorithms can be derived by using linear algebra techniques, the formal orthogonal polynomials (FOP) approach is perhaps the most common for deriving them. In fact, all recursive algorithms implementing the Lanczos method can be derived using the theory of FOP, \cite{93:Brezinski}.

A drawback of these algorithms is their inherent fragility due to the nonexistence
of some orthogonal polynomials. This causes them to breakdown well before convergence. To avoid these breakdowns, variants that jump over the nonexisting polynomials have been developed; they are referred to as breakdown-free algorithms, \cite{92:Brezinski,99:Brezinski,97:Grav,99:Guennouni,85:Parlett,94:Ye}. Note that it is not the
purpose of this paper to discuss this breakdown issue although it will be highlighted when necessary.

Two types of recurrence relations are needed: one for $P_{k}(x)$ and one for $\textit{P}^{(1)}_{k}(x)$, \cite{95:Baheux}. In \cite{94:Baheux,95:Baheux}, recurrence relations for the computation of polynomials $\textit{P}_{k}(x)$ are represented by $A_{i}$ and those for polynomials $\textit{P}^{(1)}_{k}(x)$, by $B_{j}$. Table $1$ and Table $2$ below give a comprehensive list.

C. Baheux and C. Brezinski have exploited some of the polynomial relations which involve few matrix-vector multiplications. In their work, the only relations that they studied were those where the degrees of the polynomials in the right and left hand sides of the relation differ by ONE or TWO at most. We are studying relations where the difference in degrees is TWO or THREE. For full details of these relations, see \cite{10:Farooq}. The following notation has been introduced in \cite{94:Baheux,95:Baheux}. We will adopt it here and extend the list accordingly.
\newpage
\begin{table}[h]
\caption{Computation of $A_{i}$ and $B_{j}$ from different polynomials \cite{94:Baheux}} 
\centering 
\begin{tabular}{|c|c|c|c|}
\hline
Relation $A_{i}$ & \multicolumn{1}{|l|}{Computation of $P_{k}$ from} & Relation $B_{j}$ & \multicolumn{1}{l|}{Computation of $P^{(1)}_{k}$ from}\\[0.5ex]
\hline\hline
$A_{1}$ & $P_{k-2}$ \qquad $P^{(1)}_{k-2}$  & $B_{1}$ &  $P_{k-2}$ \qquad $P^{(1)}_{k-2}$\\ 
$A_{2}$ & $P_{k-2}$ \qquad $P^{(1)}_{k-1}$  & $B_{2}$ & $P_{k-2}$ \qquad $P^{(1)}_{k-1}$\\
$A_{3}$ & $P_{k-2}$ \qquad $P^{(1)}_{k}$  & $B_{3}$ & $P_{k-2}$ \qquad $P_{k}$\\
$A_{4}$ & $P_{k-2}$ \qquad $P_{k-1}$ & $B_{4}$ & $P_{k-2}$ \qquad $P_{k-1}$\\
$A_{5}$ & $P^{(1)}_{k-2}$ \qquad $P_{k-1}$ & $B_{5}$ & $P^{(1)}_{k-2}$ \qquad $P^{(1)}_{k-1}$\\
$A_{6}$ & $P^{(1)}_{k-2}$ \qquad $P^{(1)}_{k-1}$ & $B_{6}$ & $P^{(1)}_{k-2}$ \qquad $P^{(1)}_{k-1}$\\
$A_{7}$ & $P^{(1)}_{k-2}$ \qquad $P^{(1)}_{k}$  & $B_{7}$ & $P^{(1)}_{k-2}$ \qquad $P_{k}$\\
$A_{8}$ & $P^{(1)}_{k-1}$ \qquad $P_{k-1}$ & $B_{8}$ & $P^{(1)}_{k-1}$ \qquad $P_{k-1}$\\
$A_{9}$ & $P_{k-1}$ \qquad $P^{(1)}_{k}$ & $B_{9}$ & $P_{k-1}$ \qquad $P_{k}$\\
$A_{10}$ & $P^{(1)}_{k-1}$ \qquad $P^{(1)}_{k}$  & $B_{10}$ & $P^{(1)}_{k-1}$ \qquad $P_{k}$ \\
\hline
\end{tabular}
\end{table}
\begin{table}[h]
\caption{Polynomials used in the computation of new relations $A_{i}$ and $B_{j}$} 
\centering 
\begin{tabular}{|c|c|c|c|}
\hline
Relation $A_{i}$ & \multicolumn{1}{|l|}{Computation of $P_{k}$ from} & Relation $B_{j}$ & \multicolumn{1}{l|}{Computation of $P^{(1)}_{k}$ from}\\[0.5ex]
\hline\hline
$A_{11}$ & $P_{k-3}$ \qquad $P^{(1)}_{k-1}$   & $B_{11}$ & $P_{k-3}$ \qquad $P_{k-1}$ \\ 
$A_{12}$ & $P_{k-2}$ \qquad $P_{k-3}$  &  $B_{12}$ & $P_{k-2}$ \qquad $P_{k-3}$\\
$A_{13}$ & $P_{k-2}$ \qquad $P^{(1)}_{k-3}$  &  $B_{13}$ & $P^{(1)}_{k-2}$ \qquad $P^{(1)}_{k-3}$\\
$A_{14}$ & $P^{(1)}_{k-2}$ \qquad $P^{(1)}_{k-3}$ &  $B_{14}$ & $P^{(1)}_{k-3}$ \qquad $P_{k-1}$ \\
$A_{15}$ & $P^{(1)}_{k-3}$ \qquad $P^{(1)}_{k-1}$ &  $B_{15}$ & $P_{k-2}$ \qquad $P^{(1)}_{k-2}$\\
$A_{16}$ & $P_{k-2}$ \qquad $P^{(1)}_{k-2}$ &  $B_{16}$  & $P^{(1)}_{k-2}$ \qquad $P_{k-1}$\\
$A_{17}$ & $P_{k-2}$ \qquad $P^{(1)}_{k-1}$ & -- & --\\
$A_{18}$ & $P^{(1)}_{k-1}$ \qquad $P^{(1)}_{k-2}$ & -- & --\\
$A_{19}$ & $P^{(1)}_{k-2}$ \qquad $P_{k-1}$ & -- & --\\
\hline
\end{tabular}
\end{table}

The paper is organized as follow. Section 2, briefly recalls the Lanczos algorithm. Section $3$, derives some of the possible recurrence relations $A_{i}$ and $B_{j}$ given in Table $2$ and their combination to obtain Lanczos-type algorithms. It also discusses two recurrence relations which although exist and satisfy the normalization and orthogonality conditions, cannot be used for the computation of $\textbf{r}_{k}$ and hence $\textbf{x}_{k}$. Section $4$ is the conclusion.

\section{The Lanczos algorithm}
Consider a linear system in $\textit{R}^n$ with $n$ unknowns
\begin{equation}
\textit{A}\textbf{x}=\textbf{b}.
\end{equation}
Lanczos method, \cite{95:Baheux,99:Brezinski,00:Brezinski,52:Lanczos}, for solving (1), consists in constructing a sequence of vectors $\textbf{x}_{k}$ as follows.

$\bullet$ choose two arbitrary vectors $\textbf{x}_{0}$ and $\textbf{y}\neq0$ in $\textit{R}^n$,

$\bullet$ set $\textbf{r}_{0}=\textbf{b}-\textit{A}\textbf{x}_{0}$,

$\bullet$ determine $\textbf{x}_{k}$ such that

$\textbf{x}_{k}-\textrm{x}_{0}\in\textit{E}_{k}=span(\textbf{r}_{0},\textit{A}\textbf{r}_{0},\dots,\textit{A}^{k-1}\textbf{r}_{0})$

$\textbf{r}_{k}=\textbf{b}-\textit{A}\textbf{x}_{k}\bot\textit{F}_{k}=span(\textbf{y},\textit{A}^T\textbf{y},\dots,\textit{A}^{T^{k-1}}\textbf{y})$

\noindent where $A^T$ is transpose of $A$.

\noindent $\textbf{x}_k-\textbf{x}_0$ can be written as
\[\textbf{x}_{k}-\textbf{x}_{0}=-\alpha_{1}\textbf{r}_{0}- \dots -\alpha_{k}\textit{A}^{k-1}\textbf{r}_{0}.\]
Multiplying both sides by $A$, adding and subtracting b and simplifying, we get
\[\textbf{r}_{k}=\textbf{r}_{0}+\alpha_{1}\textit{A}\textbf{r}_{0}+\dots+\alpha_{k}\textit{A}^{k}\textbf{r}_{0}\]
\noindent and the orthogonality condition above give

\[(\textit{A}^{T^{i}}\textbf{y},\textbf{r}_{k})=0 \mbox{ for i = 0,\dots, k-1},\]

\noindent which is a system of $k$ linear equations in the $k$ unknowns $\alpha_{1},\dots,\alpha_{k}$. This system is nonsingular only if $\textbf{r}_{0},\textit{A}\textbf{r}_{0},\dots,\textit{A}^{k-1}\textbf{r}_{0}$ and $\textbf{y},\textit{A}^T\textbf{y},\dots,\textit{A}^{T^{k-1}}\textbf{y}$ are linearly independent.

If we set
\[\textit{P}_{k}(\xi)=1+\alpha_{1}\xi+\dots+\alpha_{k}\xi^{k}\]
\noindent then we have \[\textbf{r}_{k}=\textit{P}_{k}(\textit{A})\textbf{r}_{0}.\]

\noindent Moreover, if we set
\[c_{i}=(\textbf{y},\textit{A}^i\textbf{r}_{0}) \mbox{ for i=0,1,\dots }\]
\noindent and we define the linear functional $c$ on the space of polynomials by
\[c(\xi^i)=c_{i} \mbox{ for i=0,1,\dots }\]
\noindent then the preceding orthogonality conditions can be written as
\[c(\xi^i\textit{P}_{k})=0 \mbox{ for i=0,\dots },k-1\]

These relations show that $\textit{P}_{k}$ is the polynomial of degree at most $k$ belonging to the family of orthogonal polynomials with respect to c, \cite{80:Brezinski}. This polynomial is defined apart from a multiplying factor which was chosen, in our case, such that $\textit{P}_{k}(0)=1$. With this normalization condition, $\textit{P}_{k}$ exists and is unique if and only if the following Hankel determinant
\[\textbf{H}^{(1)}_{k}=\left\vert\begin{array}{cccc}
c_{1} & c_{2} & \cdots & c_{k}\\
c_{2} & c_{3} & \cdots & c_{k+1}\\
\vdots & \vdots &  & \vdots\\
c_{k} & c_{k+1} & \cdots & c_{2k-1}\\
\end{array}\right\vert\]\noindent is different from zero.

Let us now consider the monic polynomial $\textit{P}^{(1)}_{k}$ of degree $k$ belonging to the FOP with respect to the functional $c^{(1)}$ defined by $c^{(1)}(\xi^i)=c(\xi^{i+1})$. $\textit{P}^{(1)}_{k}$ exists and is unique, if and only if the Hankel determinant $\textbf{H}^{(1)}_{k}\neq0$, which is the same condition as for the existence and uniqueness of $\textit{P}_{k}$.

A Lanczos-type algorithm consists in computing $\textit{P}_{k}$ recursively, then $\textbf{r}_{k}$ and finally $\textbf{x}_{k}$ such that $\textbf{r}_{k}=\textbf{b}-\textit{A}\textbf{x}_{k}$, without inverting $A$, which gives the solution of the system $(1)$ in at most $n$ steps, in exact arithmetic where $n$ is the dimension of the system of linear equations, \cite{93:Brezinski,99:Brezinski}.

\section{Recursive computation of $\textit{P}_{k}$ and $\textit{P}^{(1)}_{k}$}
The recursive computation of the polynomials $\textit{P}_{k}$, needed in the Lanczos method, can be achieved in many ways. We can use the usual three-term recurrence relation, or the relation involving the polynomials of the form $\textit{P}^{(1)}_{k}$. Such recurrence relations lead to all known algorithms for implementing Lanczos-type algorithms. For a unified presentation of all these methods based on the theory of FOP, see \cite{95:Baheux,99:Brezinski,10:Farooq}.

We need two recurrence relations, one for $\textit{P}_{k}$ and one for $\textit{P}^{(1)}_{k}$. All Lanczos-type algorithms are characterized by the choice of these recurrence relationships. In the following we will discuss some of these recurrence relations for $\textit{P}_{k}$, \cite{10:Farooq}, and derive new recurrence relationships between adjacent orthogonal polynomials,
\cite{94:Baheux,95:Baheux,80:Brezinski,93:Brezinski,92:Zaglia,00:Brezinski,39:Szego}, which can be used to design new Lanczos-type algorithms, as has been shown in \cite{95:Baheux,99:Brezinski,10:Farooq}. Note that the term ``Lanczos process'' and ``Lanczos-type algorithm'' are used interchangeably throughout the paper.

\subsection{Relations $A_{i}$}
We will follow the notation explained in Section $1$. First we will derive relations $A_{i}$ ( $i>10$ ) for $P_{k}$ which can be used to find $\textbf{r}_{k}$ and then $\textbf{x}_{k}$ without inverting $A$, the matrix of the system to be solved. We will only try to find the constant coefficients of recurrence relations which can be used for the implementation of Lanczos-type algorithms. If a recurrence relation exists but cannot be used for such an implementation, then there is no need to calculate its coefficients. The reason for that will be given. Note, however, that when a recurrence relation exists and can be computed, leading therefore to a Lanczos-type algorithm, the algorithm may still break down for two reasons:\\
 \noindent 1. There is a loss of orthogonality as in most known Lanczos-type algorithms, \cite{06:Gérard,03:Saad}.\\
 \noindent 2. When coefficients involve determinants in their denominators, these determinants may be zero (or rather, in practice, just close to zero). In this case, a check on the value of the determinants will determine whether the process has to be stopped or not. This check appears in the relevant algorithms given in \cite{10:Farooq}. This kind of breakdown is called ghost breakdown, \cite{94:Brezinski,97:Brezinski}. It may be cured by conditioning methods, \cite{55:Riley,96:Kim,00:Xue}. This is not considered here. However, it can also be cured using the restarting and switching strategies, presented in \cite{10:Farooq,11:Farooq,12:Farooq}.

 In the following $P_{k}$ stands for $P_{k}(x)$ and $P^{(1)}_{k}$ for $P^{(1)}_{k}(x)$. But, before we derive the relations $A_i$ and $B_j$, we first define the notion of``orthogonal polynomials sequence", \cite{84:Chihara}.\\
\textbf{Definition 1.}
A sequence $\{P_n\}$ is called an orthogonal polynomial sequence with respect to the linear functional $c$ if, for all nonnegative integers $m$ and $n$,\\
\noindent $(i)$ $P_n$ is polynomial of degree $n$,\\
\noindent $(ii)$ $c(x^mP_n)=0$, for $m\neq n$,\\
\noindent $(iii)$ $c(x^nP_n)\neq 0$.\\
\subsubsection{Relation $A_{11}$}
As explained earlier, we follow up from what is already known up to $A_{10}$, \cite{94:Baheux}. $A_{11}$ is therefore the natural follow up. Consider the following recurrence relationship
\begin{equation}
P_{k}(x)=(A_{k}x^{3}+B_{k}x^2+C_{k}x
+D_{k})P_{k-3}(x)+(E_{k}x+F_{k})P^{(1)}_{k-1}(x),
\end{equation}
\noindent where $P_{k}$, $P^{(1)}_{k-1}$ and $P_{k-3}$ are polynomials of degree $k$, $k-1$ and $k-3$ respectively.

\noindent {\bf Proposition 3.1:} Relation of the form $A_{11}$ does not exist.\\
\noindent {\bf Proof:}
Let us see if all coefficients can be identified. If $x^{i}$ is a polynomial of exact degree $i$ then
\begin{center}$\forall i=0, \dots ,k-1$, $c(x^{i}P_{k})=0. \longrightarrow (C_{1})$\end{center}
\begin{center}$\forall i=0, \dots ,k-2$,
$c^{(1)}(x^{i}P^{(1)}_{k-1})=0. \longrightarrow
(C_{2})$\end{center}
\begin{center}$\forall i=0, \dots ,k-4$, $c(x^{i}P_{k-3})=0. \longrightarrow (C_{3})$\end{center}
\noindent where $c$ and $c^{(1)}$ are defined respectively as follows. \begin{center}$c(x^i)=c_{i}. \longrightarrow
(C_{4})$\end{center}
\begin{center}$c^{(1)}(x^i)=c(x^{i+1}). \longrightarrow (C_{5})$\end{center}
\noindent Since $P_k(0)=1$, $\forall k$, then for $x=0$, equation (2) becomes
\begin{equation}
1=D_{k}+F_{k}P^{(1)}_{k-1}(0).
\end{equation}

\noindent Multiplying both sides of equation (2) by $x^{i}$ and applying $c$ using the condition $(C_{5})$ where necessary, it can be written as
\begin{align}
&c(x^{i}P_{k})=A_{k}c(x^{i+3}P_{k-3})+B_{k}c(x^{i+2}P_{k-3}) +C_{k}c(x^{i+1}P_{k-3})\\ \nonumber
&               \quad \quad\qquad+D_{k}c(x^{i}P_{k-3}) +E_{k}c^{(1)}(x^{i}P^{(1)}_{k-1})+F_{k}c(x^{i}P^{(1)}_{k-1}).
\end{align}
\noindent For $i=0$, equation (4) gives $0=F_{k}c(P^{(1)}_{k-1})$. Since $c(P^{(1)}_{k-1})\neq0$, \cite{80:Brezinski,91:Brezinski} then $F_{k}=0$. Hence from (3), we have $D_{k}=1$.

\noindent Equation (4) is always true for $i=1,...,k-7$ by $(C_{2})$ and $(C_{3})$.

\noindent For $i=k-6$, (4) becomes $A_{k}c(x^{k-3}P_{k-3})=0$, as all other terms vanish due to conditions $(C_{1})$, $(C_{2})$ and $(C_{3})$. But, according to condition $(iii)$ of definition $1$ in Section $3.1$, $c(x^{k-3}P_{k-3})\neq0$, therefore, $A_{k}=0$.

\noindent For $i=k-5$, (4) becomes $B_{k}c(x^{k-3}P_{k-3})=0$. Since $c(x^{k-3}P_{k-3})\neq0$, $B_{k}=0$.

\noindent For $i=k-4$, (4) becomes $C_{k}c(x^{k-3}P_{k-3})=0$. Since $c(x^{k-3}P_{k-3})\neq0$, $C_{k}=0$.

\noindent Putting values of $A_{k}$, $B_{k}$, $C_{k}$, $D_{k}$ and $F_{k}$ in equation (4), we get
\[c(x^{i}P_{k})=c(x^{i}P_{k-3})
+E_{k}c^{(1)}(x^{i}P^{(1)}_{k-1}).\]

\noindent For $i=k-3$, this equation becomes

\[c(x^{k-3}P_{k})=c(x^{k-3}P_{k-3})
+E_{k}c^{(1)}(x^{k-3}P^{(1)}_{k-1}).\]

\noindent Using $(C_{1})$ and $(C_{2})$, both $c(x^{k-3}P_{k})=0$ and $c^{(1)}(x^{k-3}P^{(1)}_{k-1})=0$, therefore we get
\[c(x^{k-3}P_{k-3})=0,\]
\noindent which is impossible due to condition $(iii)$ of definition 1 given in Section $3.1$. Therefore, Proposition $3.1$ holds.

\subsubsection{Relation $A_{12}$}
Consider the following recurrence relationship for $k\geq3$,
\begin{equation}
P_{k}(x)=A_{k}[(x^{2}+B_{k}x+C_{k})P_{k-2}(x)+(D_{k}x^3+E_{k}x^2+F_{k}x+G_{k})P_{k-3}(x)],
\end{equation}

This recurrence relation has been considered in \cite{10:Farooq,10:Salhi}.
\subsubsection{Relation $A_{13}$}
Consider the following recurrence relationship
\begin{equation}\label{A13}
P_{k}(x)=A_{k}[(x^{2}+B_{k}x+C_{k})P_{k-2}(x)
+(D_{k}x^3+E_{k}x^2+F_{k}x+G_{k})P^{(1)}_{k-3}(x)],
\end{equation}
\noindent where $P_{k}$, $P_{k-2}$ and $P^{(1)}_{k-3}$ are orthogonal polynomials of degree $k$, $k-2$ and $k-3$ respectively and $A_{k}$, $B_{k}$, $C_{k}$, $D_{k}$, $E_{k}$, $F_{k}$ and $G_{k}$ are constants to be determined using the normalization condition $P_{k}(0)=1$ and the orthogonality condition $(C_{1})$.

\noindent{\bf Proposition 3.2:} Relation of the form $A_{13}$ exists.

\noindent{\bf Proof:}
We know that
\begin{center}$\forall  i=0, \dots ,k-1$, $c^{(1)}(x^{i}P^{(1)}_{k})=0$. $\longrightarrow$ $(C_{6})$\end{center}
Since $\forall k, P_{k}(0)=1$, equation \eqref{A13} gives
\begin{equation}\label{7}
1=A_{k}[C_{k}+G_{k}P^{(1)}_{k-3}(0)].
\end{equation}
\noindent Multiplying both sides of \eqref{A13} by $x^{i}$ and then applying the linear functional $c$, we get
\begin{equation}\begin{array}{l}\label{8}
c(x^{i}P_{k})=A_{k}[c(x^{i+2}P_{k-2})+B_{k}c(x^{i+1}P_{k-2})+C_{k}c(x^{i}P_{k-2})\\+D_{k}c^{(1)}(x^{i+2}P^{(1)}_{k-3})+E_{k}c^{(1)}(x^{i+1}P^{(1)}_{k-3})+F_{k}c^{(1)}(x^{i}P^{(1)}_{k-3})+G_{k}c(x^{i}P^{(1)}_{k-3})].
\end{array}\end{equation}
\noindent For $i=0$, equation \eqref{8} becomes, $0=G_{k}c(P^{(1)}_{k-3})$. Since $c(P^{(1)}_{k-3})\neq0$, this implies that $G_{k}=0.$ Therefore, from equation \eqref{7} we have $A_{k}=\frac{1}{C_{k}}$.

\noindent The orthogonality condition $(C_1)$ is always true for $i=1,\dots,k-6$.

\noindent For $i=k-5$, we have $D_{k}c^{(1)}(x^{k-3}P^{(1)}_{k-3})=0$, which implies that $D_{k}=0$.

\noindent For $i=k-4$, \eqref{8} becomes $c(x^{k-2}P_{k-2})+E_{k}c^{(1)}(x^{k-3}P^{(1)}_{k-3})=0$,
\[E_{k}=-\frac{c(x^{k-2}P_{k-2})}{c^{(1)}(x^{k-3}P^{(1)}_{k-3})}.\]
\noindent For $i=k-3$, \eqref{8} gives
\begin{equation}\begin{array}{l}\label{9}
B_{k}c(x^{k-2}P_{k-2})+F_{k}c^{(1)}(x^{k-3}P^{(1)}_{k-3})\\=-c(x^{k-1}P_{k-2})
-E_{k}c^{(1)}(x^{k-2}P^{(1)}_{k-3}).
\end{array}\end{equation}
\noindent For $i=k-2$, \eqref{8} becomes
\begin{equation}\begin{array}{l}\label{10}
B_{k}c(x^{k-1}P_{k-2})+C_{k}c(x^{k-2}P_{k-2})+F_{k}c^{(1)}(x^{k-2}P^{(1)}_{k-3})
\\=-c(x^{k}P_{k-2})-E_{k}c^{(1)}(x^{k-1}P^{(1)}_{k-3}).
\end{array}\end{equation}
\noindent For $i=k-1$, we get
\begin{equation}\begin{array}{l}\label{11}
B_{k}c(x^{k}P_{k-2})+C_{k}c(x^{k-1}P_{k-2})+F_{k}c^{(1)}(x^{k-1}P^{(1)}_{k-3})\\
=-c(x^{k+1}P_{k-2})-E_{k}c^{(1)}(x^{k}P^{(1)}_{k-3}).
\end{array}\end{equation}
Let $a_{11}$, $a_{21}$, $a_{31}$, $a_{12}$, $a_{22}$, $a_{32}$ and $a_{13}$, $a_{23}$, $a_{33}$ be the coefficients of $B_{k}$, $C_{k}$ and $G_{k}$ in equations \eqref{9}, \eqref{10} and \eqref{11} respectively and suppose $b_{1}$, $b_{2}$ and $b_{3}$ are the corresponding right hand side terms of these equations. Then we have

\noindent $a_{11}=c(x^{k-2}P_{k-2})$, $a_{12}=0$, $a_{13}=c^{(1)}(x^{k-3}P^{(1)}_{k-3})$,

\noindent $a_{21}=c(x^{k-1}P_{k-2})$, $a_{22}=c(x^{k-2}P_{k-2})$, $a_{23}=c^{(1)}(x^{k-2}P^{(1)}_{k-3})$,

\noindent $a_{31}=c(x^{k}P_{k-2})$, $a_{32}=c(x^{k-1}P_{k-2})$, $a_{33}=c^{(1)}(x^{k-1}P^{(1)}_{k-3})$,

\noindent
$b_{1}=-c(x^{k-1}P_{k-2})-E_{k}c^{(1)}(x^{k-2}P^{(1)}_{k-3})$,
$b_{2}=-c(x^{k}P_{k-2})-E_{k}c^{(1)}(x^{k-1}P^{(1)}_{k-3})$,

\noindent $b_{3}=-c(x^{k+1}P_{k-2})-E_{k}c^{(1)}(x^{k}P^{(1)}_{k-3})$,
\begin{equation}
a_{11}B_{k}+0C_{k}+a_{13}F_{k}=b_{1},
\end{equation}
\begin{equation}
a_{21}B_{k}+a_{22}C_{k}+a_{23}F_{k}=b_{2},
\end{equation}
\begin{equation}
a_{31}B_{k}+a_{32}C_{k}+a_{33}F_{k}=b_{3}.
\end{equation}
If $\Delta_{k}$ represents the determinant of the coefficients matrix of the above system of equations then we have
\[\Delta_{k}=a_{11}(a_{22}a_{33}-a_{32}a_{23})+a_{13}(a_{21}a_{32}-a_{31}a_{22}).\]
\noindent If $\Delta_{k}\neq0$, then
\[B_{k}=\frac{[b_{1}(a_{22}a_{33}-a_{32}a_{23})+a_{13}(b_{2}a_{32}-b_{3}a_{22})]}{\Delta_{k}},\]
\[F_{k}=\frac{b_{1}-a_{11}B_{k}}{a_{13}},\]
\[C_{k}=\frac{b_{2}-a_{21}B_{k}-a_{23}F_{k}}{a_{22}}\] and
\[A_{k}=\frac{1}{C_{k}}.\]
\noindent Hence,
\begin{equation}P_{k}(x)=A_{k}[(x^{2}+B_{k}x+C_{k})P_{k-2}(x)+(E_{k}x^2+F_{k}x)P^{(1)}_{k-3}(x)].\end{equation}
Here again, the $A_{13}$ recurrence relation exists. It therefore, can be used to implement the Lanczos process.\\

\noindent {\bf Remark 1:} If $\Delta_k = 0$ then we cannot estimate coefficient $B_k$, which means
relation $A_{13}$ may not exit. However, it may exist, but for one or more of its
coefficients which cannot be estimated for numerical reasons. This would be
a case of ghost breakdown.

\subsubsection{Relation $A_{14}$}
Consider the following recurrence relationship
\begin{equation}\label{A14}
P_{k}(x)=A_{k}[(x^{2}+B_{k}x+C_{k})P^{(1)}_{k-2}(x)
+(D_{k}x^3+E_{k}x^2+F_{k}x+G_{k})P^{(1)}_{k-3}(x)],
\end{equation}
\noindent where $P_{k}$, $P^{(1)}_{k-2}$ and $P^{(1)}_{k-3}$ are orthogonal polynomials of degree $k$, $k-2$ and $k-3$ respectively and $A_{k}$, $B_{k}$, $C_{k}$, $D_{k}$, $E_{k}$, $F_{k}$ and $G_{k}$ are constants to be determined using the normalization condition $P_{k}(0)=1$ and the orthogonality conditions $(C_{1})$ and $(C_{6})$.

\noindent{\bf Proposition 3.3:} Relation of the form $A_{14}$ exists.

\noindent{\bf Proof:}
Since $\forall k, P_{k}(0)=1$, equation \eqref{A14} gives
\begin{equation}\label{17}
1=A_{k}[C_{k}P^{(1)}_{k-2}(0)+G_{k}P^{(1)}_{k-3}(0)].
\end{equation}
\noindent Multiplying both sides of equation \eqref{A14} by $x^{i}$ and then applying the linear functional $c$ and using condition
$(C_{5})$, we get
\begin{equation}\begin{array}{l}\label{18}
c(x^{i}P_{k})=A_{k}[c^{(1)}(x^{i+1}P^{(1)}_{k-2})+B_{k}c^{(1)}(x^{i}P^{(1)}_{k-2})
+C_{k}c(x^{i}P^{(1)}_{k-2})\\+D_{k}c^{(1)}(x^{i+2}P^{(1)}_{k-3})+E_{k}c^{(1)}(x^{i+1}P^{(1)}_{k-3})
+F_{k}c^{(1)}(x^{i}P^{(1)}_{k-3})+G_{k}c(x^{i}P^{(1)}_{k-3})].
\end{array}\end{equation}
\noindent For $i=0$, equation \eqref{18} becomes
\begin{equation}\label{19}
C_{k}c(P^{(1)}_{k-2})+G_{k}c(P^{(1)}_{k-3})=0.
\end{equation}
\noindent The orthogonality condition for \eqref{18} is always true for $i=1, \dots ,k-6$.

\noindent For $i=k-5$, we get $D_{k}c^{(1)}(x^{k-3}P^{(1)}_{k-3})=0$. Since $c^{(1)}(x^{k-3}P^{(1)}_{k-3})\neq0$, $D_{k}=0$.

\noindent For $i=k-4$, $E_{k}c^{(1)}(x^{k-3}P^{(1)}_{k-3})=0$. Since $c^{(1)}(x^{k-3}P^{(1)}_{k-3})\neq0$, we have $E_{k}=0$.

\noindent For $i=k-3$, \begin{equation}\label{20}
F_{k}c^{(1)}(x^{k-3}P^{(1)}_{k-3})+G_{k}c(x^{k-3}P^{(1)}_{k-3})=-c^{(1)}(x^{k-2}P^{(1)}_{k-2}). \end{equation}
\noindent For $i=k-2$, \begin{equation}\begin{array}{l}\label{21}
B_{k}c^{(1)}(x^{k-2}P^{(1)}_{k-2})+C_{k}c(x^{k-2}P^{(1)}_{k-2})
+F_{k}c^{(1)}(x^{k-2}P^{(1)}_{k-3})\\+G_{k}c(x^{k-2}P^{(1)}_{k-3})
=-c^{(1)}(x^{k-1}P^{(1)}_{k-2}). \end{array}\end{equation}
\noindent For $i=k-1$, \begin{equation}\begin{array}{l}\label{22}
B_{k}c^{(1)}(x^{k-1}P^{(1)}_{k-2})+C_{k}c(x^{k-1}P^{(1)}_{k-2})+
F_{k}c^{(1)}(x^{k-1}P^{(1)}_{k-3})\\+G_{k}c(x^{k-1}P^{(1)}_{k-3})
=-c^{(1)}(x^{k}P^{(1)}_{k-2}). \end{array}\end{equation} \noindent The values of $A_{k}$, $B_{k}$, $C_{k}$, $F_{k}$ and $G_{k}$ can be obtained by solving equations \eqref{17}, \eqref{19}, \eqref{20}, \eqref{21} and \eqref{22}. Hence
\[P_{k}(x)=A_{k}[(x^{2}+B_{k}x+C_{k})P^{(1)}_{k-2}(x)+(F_{k}x+G_{k})P^{(1)}_{k-3}(x)].\]
Now multiplying both sides of the above relation by $\textbf{r}_{0}$, replacing $x$ by $A$ and using the relations $\textbf{r}_{k}=P_{k}(A)\textbf{r}_{0}$ and $\textbf{z}_{k}=P^{(1)}_{k}(A)\textbf{r}_{0}$, we get
\begin{equation}
\textbf{r}_{k}=A_{k}[(\textit{A}^{2}+B_{k}\textit{A}+C_{k}\textit{I})\textbf{z}_{k-2}
+(F_{k}\textit{A}+G_{k}\textit{I})\textbf{z}_{k-3}].
\end{equation}
\noindent Using $\textbf{r}_{k}=\textbf{b}-\textit{A}\textbf{x}_{k}$, we get from the last equation
\begin{equation}
\textit{A}\textbf{x}_{k}=\textbf{b}-A_{k}[(\textit{A}^{2}+B_{k}\textit{A}
+C_{k}\textit{I})\textbf{z}_{k-2}+(F_{k}\textit{A}+G_{k}\textit{I})\textbf{z}_{k-3}].
\end{equation}
From this relation it is clear that we cannot find $\textbf{x}_{k}$ from $\textbf{r}_{k}$ without inverting
A. Hence, this recurrence relation exist as stipulated by Proposition 3.3. But, it is not desirable to implement a Lanczos-type algorithm.

For recurrence relations $A_{15}$, $A_{16}$, $A_{17}$, $A_{18}$ and $A_{19}$ and their corresponding
coefficients, consult \cite{10:Farooq}.

\subsection{Relations $B_{j}$}
Now we consider relations of the type $B_{j}$ which have not been considered before, \cite{94:Baheux,95:Baheux}, i.e $B_{j}$ with $j>10$. These relations, when they exist will be used in combination with relations $A_{i}$ to derive further Lanczos-type algorithms as explained in \cite{10:Farooq}.

\subsubsection{Relation $B_{11}$}
Consider the following recurrence relationship
\begin{equation}\label{B11}
P^{(1)}_{k}(x)=(A_{k}x^3+B_{k}x^2+
C_{k}x+D_{k})P_{k-3}(x)+(E_{k}x+F_{k})P_{k-1}(x),\end{equation}
\noindent where $P^{(1)}_{k}(x)$, $P_{k-1}(x)$ and $P_{k-3}(x)$ are orthogonal polynomials of degree $k$, $k-1$ and $k-3$ respectively.

\noindent{\bf Proposition 3.4:} Relation of the form $B_{11}$ does not exist.

\noindent{\bf Proof:}
Let $x^{i}$ be a polynomial of exact degree $i$ then \begin{center}$\forall i=0, \dots ,k-4$, $c(x^{i}P_{k-3})=0. \longrightarrow (C_{10})$\end{center}
Multiply both sides of equation \eqref{B11} by $x^{i}$ and applying $c^{(1)}$ and also using condition $(C_{5})$ where necessary, we get
\begin{equation}\begin{array}{l}\label{26}
c^{(1)}(x^{i}P^{(1)}_{k})=A_{k}c(x^{i+4}P_{k-3})+B_{k}c(x^{i+3}P_{k-3})+C_{k}c(x^{i+2}P_{k-3})\\
+D_{k}c(x^{i+1}P_{k-3})+E_{k}c(x^{i+2}P_{k-1})+F_{k}c(x^{i+1}P_{k-1}).
\end{array}\end{equation}
\noindent The relation \eqref{26} is always true for $i=0,\dots,k-8.$

\noindent For $i=k-7$, we have $A_{k}c(x^{k-3}P_{k-3})=0$. Which implies $A_{k}=0$, because $c(x^{k-3}P_{k-3})\neq0.$

Similarly for $i=k-6$, $i=k-5$, $i=k-4$ and $i=k-3$, we get respectively $B_{k}=0$, $C_{k}=0$, $D_{k}=0$ and $E_{k}=0$. But $P^{(1)}_{k}(x)$ is a monic polynomial of degree $k$ and we see that $A_{k}=0$. Therefore $E_{k}=1$. If $E_{k}=0$ then $P^{(1)}_{k}(x)$ is no more of degree $k$ as the degree of the $P^{(1)}_{k}$ depends on $A_{k}$ and $E_{k}$ and we know that $A_{k}=0$ and if $E_{k}=1$ then $c(x^{k-1}P_{k-1})=0$ which is also impossible. Similarly for $i=k-2$, we get $F_{k}=0$ and if $E_{k}=0$, then $P^{(1)}_{k}=0.$ Hence the relation $B_{11}$ does not exist, therefore, Proposition $3.4$ holds.

\medskip
Recurrence relation $B_{12}$ has been explored in \cite{10:Farooq}.

\subsubsection{Relation $B_{13}$}
Consider the following recurrence relationship
\begin{equation}\label{B13}
P^{(1)}_{k}(x)=(A_{k}x^3+B_{k}x^2+C_{k}x+
D_{k})P^{(1)}_{k-3}(x)+(E_{k}x^2+F_{k}x+G_{k})P^{(1)}_{k-2}(x),
\end{equation}
\noindent where $P^{(1)}_{k}(x)$, $P^{(1)}_{k-2}(x)$ and $P^{(1)}_{k-3}(x)$ are orthogonal polynomials of degree $k$, $k-2$ and $k-3$ respectively.

\noindent{\bf Proposition 3.5:} Relation of the form $B_{13}$ exists.

\noindent{\bf Proof:}
Let $x^{i}$ be a polynomial of exact degree $i$ then
\begin{center}$\forall i=0, \dots ,k-3$,
$c^{(1)}(x^{i}P^{(1)}_{k-2})=0. \longrightarrow
(C_{11})$\end{center}
Multiply both sides of equation \eqref{B13} by $x^{i}$ and applying $c^{(1)}$, we get
\begin{equation}\begin{array}{l}\label{28}
c^{(1)}(x^{i}P^{(1)}_{k})=A_{k}c^{(1)}(x^{i+3}P^{(1)}_{k-3})
+B_{k}c^{(1)}(x^{i+2}P^{(1)}_{k-3})+C_{k}c^{(1)}(x^{i+1}P^{(1)}_{k-3})
\\+D_{k}c^{(1)}(x^{i}P^{(1)}_{k-3})+E_{k}c^{(1)}(x^{i+2}P^{(1)}_{k-2})
+F_{k}c^{(1)}(x^{i+1}P^{(1)}_{k-2})+G_{k}c^{(1)}(x^{i}P^{(1)}_{k-2}).
\end{array}\end{equation}
\noindent The orthogonality condition is always true for $i=0,...,k-7$.

\noindent For $i=k-6$, we get $A_{k}c^{(1)}(x^{k-3}P^{(1)}_{k-3})=0$, which implies that $A_{k}=0$ as \\$c^{(1)}(x^{k-3}P^{(1)}_{k-3})\neq0.$ But $P^{(1)}_{k}(x)$ is a monic polynomial of degree $k$. Therefore $E_{k}=1$.

\noindent For $i=k-5$, we get $B_{k}c^{(1)}(x^{k-3}P^{(1)}_{k-3})=0$. Since $c^{(1)}(x^{k-3}P^{(1)}_{k-3})\neq0$, $B_{k}=0$.

\noindent For $i=k-4$, we have
\[C_{k}=-\frac{c^{(1)}(x^{k-2}P^{(1)}_{k-2})}{c^{(1)}(x^{k-3}P^{(1)}_{k-3})}.\]
\noindent For $i=k-3$, we get
\begin{equation}\label{29}
D_{k}c^{(1)}(x^{k-3}P^{(1)}_{k-3})+F_{k}c^{(1)}(x^{k-2}P^{(1)}_{k-2})
=-c^{(1)}(x^{k-1}P^{(1)}_{k-2})-C_{k}c^{(1)}(x^{k-2}P^{(1)}_{k-3}).
\end{equation}
\noindent For $i=k-2$, \eqref{28} becomes
\begin{equation}\begin{array}{r}\label{30}
D_{k}c^{(1)}(x^{k-2}P^{(1)}_{k-3})+F_{k}c^{(1)}(x^{k-1}P^{(1)}_{k-2})
+G_{k}c^{(1)}(x^{k-2}P^{(1)}_{k-2})\\=-c^{(1)}(x^{k}P^{(1)}_{k-2})
-C_{k}c^{(1)}(x^{k-1}P^{(1)}_{k-3}).
\end{array}\end{equation}
\noindent For $i=k-1$, \eqref{28} gives
\begin{equation}\begin{array}{r}\label{31}
D_{k}c^{(1)}(x^{k-1}P^{(1)}_{k-3})+F_{k}c^{(1)}(x^kP^{(1)}_{k-2})
+G_{k}c^{(1)}(x^{k-1}P^{(1)}_{k-2})\\=-c^{(1)}(x^{k+1}P^{(1)}_{k-2})
-C_{k}c^{(1)}(x^{k}P^{(1)}_{k-3}).
\end{array}\end{equation}
\noindent Let $a'_{11}=c^{(1)}(x^{k-3}P^{(1)}_{k-3})$, using $(C_{5})$, $a'_{11} =c(x^{k-2}P^{(1)}_{k-3})$. By the same condition we can write,

\noindent $a'_{12}=c^{(1)}(x^{k-2}P^{(1)}_{k-2})=c(x^{k-1}P^{(1)}_{k-2})$, $a'_{13}=0$,

\noindent
$a'_{21}=c^{(1)}(x^{k-2}P^{(1)}_{k-3})=c(x^{k-1}P^{(1)}_{k-3})$,
$a'_{22}=c^{(1)}(x^{k-1}P^{(1)}_{k-2})=c(x^{k}P^{(1)}_{k-2})$,

\noindent
$a'_{23}=c^{(1)}(x^{k-2}P^{(1)}_{k-2})=a'_{12}$,
$a'_{31}=c^{(1)}(x^{k-1}P^{(1)}_{k-3})=c(x^{k}P^{(1)}_{k-3})$,

\noindent
$a'_{32}=c^{(1)}(x^{k}P^{(1)}_{k-2})=c(x^{k+1}P^{(1)}_{k-2})$,
$a'_{33}=c^{(1)}(x^{k-1}P^{(1)}_{k-2})=a'_{22}$,

\noindent $b'_{1}=-c^{(1)}(x^{k-1}P^{(1)}_{k-2})-C_{k}c^{(1)}(x^{k-2}P^{(1)}_{k-3})=-a'_{22}-a'_{21}C_{k}$,

\noindent $b'_{2}=-c^{(1)}(x^{k}P^{(1)}_{k-2})-C_{k}c^{(1)}(x^{k-1}P^{(1)}_{k-3})=-a'_{32}-a'_{31}C_{k}$,

\noindent
$b'_{3}=-c^{(1)}(x^{k+1}P^{(1)}_{k-2})-C_{k}c^{(1)}(x^{k}P^{(1)}_{k-3})$.

\noindent Then equations \eqref{29}, \eqref{30} and \eqref{31} become
\begin{equation}\label{32}
a'_{11}D_{k}+a'_{12}F_{k}=b'_{1},
\end{equation}
\begin{equation}\label{33}
a'_{21}D_{k}+a'_{22}F_{k}+a'_{23}G_{k}=b'_{2}
\end{equation} \noindent and
\begin{equation}\label{34}
a'_{31}D_{k}+a'_{32}F_{k}+a'_{33}G_{k}=b'_{3}.
\end{equation}
\noindent If $\Delta'_{k}$ is the determinant of the coefficient matrix of the equations \eqref{32}, \eqref{33} and \eqref{34} then
\[\Delta'_{k}=a'_{11}(a'_{22}a'_{33}-a'_{32}a'_{23})-a'_{12}(a'_{21}a'_{33}-a'_{31}a'_{23}).\]
\noindent If $\Delta'_{k}\neq0$, then
\[D_{k}=\frac{b'_{1}(a'_{22}a'_{33}-a'_{32}a'_{23})-a'_{12}(b'_{2}a'_{33}-b'_{3}a'_{23})}{\Delta'_{k}},\]
\[F_{k}=\frac{b'_{1}-a'_{11}D_{k}}{a'_{12}}\]
\noindent and
\[G_{k}=\frac{b'_{2}-a'_{21}D_{k}-a'_{22}F_{k}}{a'_{23}}.\]
\noindent Hence, relation \eqref{B13} can be written as
\begin{equation}P^{(1)}_{k}(x)=(C_{k}x+D_{k})P^{(1)}_{k-3}(x)+(x^2+F_{k}x+G_{k})P^{(1)}_{k-2}(x),\end{equation}
\noindent and, therefore, exists as stipulated in Proposition 3.5.

\noindent {\bf Remark 2:} For the case where $\Delta'_k=0$, please consult Remark 1 above.

\medskip
\noindent For recurrence relations $B_{14}$, $B_{15}$ and $B_{16}$ and their corresponding coefficients,
see \cite{10:Farooq}.

\section{Conclusion}
In this paper, we looked in a systematic way at new recurrence relations between FOP's which have not been considered before. In particular, we have shown that relations $A_{11}$, $A_{17}$, $B_{11}$ and $B_{12}$ do not exist; relations $A_{14}$, $A_{15}$, $A_{18}$ and $B_{14}$ exist but are not suitable for implementing new Lanczostype algorithms; and relations $A_{12}$, $A_{13}$, $A_{16}$, $A_{19}$, $B_{13}$, $B_{15}$ and $B_{16}$ exist and can be used for the implementation of new Lanczos-type algorithms, \cite{10:Farooq}. Relation $A_{12}$ is self-sufficient and leads to a new Lanczos-type algorithm
on its own, \cite{10:Salhi}, while the rest of the relations can lead to Lanczos-type algorithms when combined in $A_i/B_j$ fashion. Possible combinations, which are studied in \cite{10:Farooq,10:Salhi}, are:

\begin{center}
$A_{13}/B_{13}$, $A_{13}/B_{15}$, $A_{13}/B_{16}$,

$A_{16}/B_{13}$, $A_{16}/B_{15}$, $A_{16}/B_{16}$,

$A_{19}/B_{13}$, $A_{19}/B_{15}$, $A_{19}/B_{16}$.
\end{center}

\bigskip
\bigskip
\medskip

%
%
%
%
%
%
%
%
%
%
%
%
%
\newpage
\bibliography{sampleJPRM}
\end{document}